\documentclass[12pt]{article}
\usepackage{a4wide}
\usepackage{theorem}
\usepackage{amsfonts}
\usepackage{enumerate}
\usepackage{amsmath}
\usepackage{amscd}
\usepackage{times}
\usepackage{isolatin1}
\usepackage[all]{xy}

\theoremstyle{change}
\theorembodyfont{\sl}
\newtheorem{thm}[subsection]{Theorem.}
\newtheorem{prop}[subsection]{Proposition.}
\newtheorem{lem}[subsection]{Lemma.}

\theorembodyfont{\rmfamily}
\newtheorem{defi}[subsection]{Definition.}
\newtheorem{rem}[subsection]{Remark.}

\newtheorem{question}[subsection]{Question.}



\makeatletter 
\newenvironment{subeqn}{\refstepcounter{subsubsection}
$$}{\leqno{\rm(\thesubsubsection)}$$\global\@ignoretrue}
\makeatother 





\newenvironment{prf}[1]{\trivlist
\item[\hskip\labelsep{\it
#1.\hspace*{.3em}}]}{~\hspace{\fill}~$\square$\endtrivlist}

\newenvironment{proof}{\begin{prf}{\bf Proof}}{\end{prf}}
\newcommand{\ZZ}{{\mathbb Z}}
\newcommand{\FF}{{\mathbb F}}
\newcommand{\QQ}{{\mathbb Q}}
\newcommand{\PP}{{\mathbb P}}
\newcommand{\RR}{{\mathbb R}}
\newcommand{\CC}{{\mathbb C}}
\renewcommand{\AA}{{\mathbb A}}
\newcommand{\AAf}[1][]{\AA_{{#1}{\rm f}}}
\renewcommand{\SS}{{\mathbb S}}
\newcommand{\HH}{{\mathbb H}}
\newcommand{\HHpm}{\HH^{\pm}}
\newcommand{\Pic}{{\rm Pic}}
\newcommand{\eps}{\varepsilon}
\newcommand{\GL}{{\rm GL}}
\newcommand{\PGL}{{\rm PGL}}
\newcommand{\SL}{{\rm SL}}
\newcommand{\Aut}{{\rm Aut}}
\newcommand{\ol}{\overline}
\newcommand{\calO}{{\cal O}}
\newcommand{\calL}{{\cal L}}
\newcommand{\calY}{{\cal Y}}
\newcommand{\End}{{\rm End}}
\newcommand{\Sym}{{\rm Sym}}
\newcommand{\Gal}{{\rm Gal}}
\newcommand{\Qbar}{{\ol{\QQ}}}
\newcommand{\discr}{{\rm discr}}
\newcommand{\Norm}{{\rm Norm}}
\newcommand{\Frob}{{\rm Frob}}
\newcommand{\Irr}{\mathrm{Irr}}
\newcommand{\Lie}{\mathrm{Lie}}
\newcommand{\SO}{\mathrm{SO}}
\newcommand{\wt}{\widetilde}

\begin{document}
\title{Special points on products of modular curves.}
\author{Bas Edixhoven\footnote{partially supported by the 
Institut Universitaire de France, and by the European TMR Network 
Contract ERB FMRX 960006 ``arithmetic algebraic geometry''.}}
\date{February 10, 2003}

\maketitle

\begin{center}
In memory of my mother, Lida Edixhoven-van Elzakker (1930--2003).
\end{center}

\begin{abstract}
We prove the Andr\'e--Oort conjecture on special points of Shimura
varieties for arbitrary products of modular curves, assuming the
Generalized Riemann Hypothesis. More explicitly, this means the
following. Let $n\geq0$, and let $\Sigma$ be a subset of~$\CC^n$
consisting of points all of whose coordinates are $j$-invariants of
elliptic curves with complex multiplications. Then we prove (under
GRH) that the irreducible components of the Zariski closure
of~$\Sigma$ are \emph{special sub-varieties}, i.e., determined by
isogeny conditions on coordinates and pairs of coordinates. A weaker
variant (Thm.~\ref{thm1.3}) is proved unconditionally.

\end{abstract}

\section{Introduction.}\label{sec1}
The main goal of this article is to prove the Andr\'e-Oort conjecture
for arbitrary products of modular curves, assuming the generalized
Riemann hypothesis (GRH) for imaginary quadratic fields. This
conjecture is usually formulated for arbitrary Shimura varieties; see
\cite{Edixhoven2} for a precise statement in the general case, and the
end of this introduction for a list of results that have been proved
so far.  The conjecture in question says that the irreducible
components of the Zariski closure of any set of special points in a
Shimura variety are sub-varieties of Hodge type. In order to be
reasonably elementary in this article, we do not use the general
formalism of Shimura varieties and their sub-varieties of Hodge type
but rather state our results in more explicit terms. In fact, we will
use the same terminology as in \cite{Edixhoven1}, which deals with the
case of products of two modular curves.  (In Section~\ref{sec2} we do
use some Shimura variety formalism, but the result in that section is
only included to show that our explicit result is in fact equivalent
to the Andr\'e-Oort conjecture.)

Let $\HH$ denote the complex upper half plane, with its
$\SL_2(\RR)$-action given by fractional linear transformations.  For
$\Gamma$ a congruence subgroup of $\SL_2(\ZZ)$, we denote by
$X_\Gamma$ the complex modular curve $\Gamma\backslash\HH$, or, more
precisely, the complex algebraic curve associated to this complex
analytic variety, and we let $\pi_\Gamma$ be the quotient map from
$\HH$ to~$X_\Gamma$. We view $X_\Gamma$ as the set of isomorphism
classes of complex elliptic curves with a level structure of
type~$\Gamma$. The endomorphism ring $\End(E)$ of a complex elliptic
curve $E$ is either $\ZZ$ or an order in an imaginary quadratic
extension of~$\QQ$; in the second case $E$ is said to be a CM elliptic
curve (CM meaning complex multiplication). A point on some $X_\Gamma$
is called a CM point if the corresponding elliptic curve has~CM. A
point on a product of curves of the form $X_\Gamma$ is called a CM
point if all its coordinates are CM points.

\begin{defi}\label{def1.1}
Let $S$ be a finite set. For every $s$ in $S$, let $\Gamma_s$ be a
congruence subgroup of $\SL_2(\ZZ)$, and let $X$ be the product of 
the~$X_{\Gamma_s}$. A closed irreducible algebraic sub-variety $Z$ of 
$X$ is called special if $S$ has a partition $(S_1,\ldots,S_r)$ such that 
$X$ is the product of sub-varieties $Z_i$ of the 
$X_i:=\prod_{s\in S_i}X_{\Gamma_s}$, each of one of the forms: 
\begin{enumerate}
\item $S_i$ is a one element set, and $Z_i$ a CM point; 
\item the image of $\HH$ in $X_i$ under the map sending $\tau$ in
$\HH$ to $\pi_{\Gamma_s}(g_s\tau)$ for every $s$ in $S_i$, with the
$g_s$ elements of $\GL_2(\QQ)$ with positive determinant. 
\end{enumerate}
\end{defi}
In Section~\ref{sec2} we will show that our ad hoc notion of ``special
sub-variety of $X$'' is the same as that of ``sub-variety of Hodge
type of~$X$''. We note that a point in $X$ as above is special if and
only if it is a CM point. We say that two points $x$ and $x'$ in $X$
are isogeneous if the corresponding products of elliptic curves are
isogeneous. (We could have asked the isogenies to preserve the product
structure, but for Theorems~\ref{thm1.2} and~\ref{thm1.3} below that
would not change anything; just note that the category of elliptic
curves up to isogeny is semi-simple, and that the symmetric
group~$S_n$ is finite.)  The main results of this article are the
following two (the second is motivated by possible applications in
transcendence theory, see~\cite{CohenWustholz1}, and by work of Vatsal
and Cornut, see \cite{Cornut1} and also \cite{Edixhoven3}).
\begin{thm}\label{thm1.2}
Let $\Sigma$ be a set of special points in a finite product of modular
curves. Assume GRH for imaginary quadratic fields. Then all
irreducible components of the Zariski closure of $\Sigma$ are special.
\end{thm}
\begin{thm}\label{thm1.3}
Let $\Sigma$ be a set of special points in a finite product of modular
curves, lying in one isogeny class.  Then all irreducible components
of the Zariski closure of $\Sigma$ are special.
\end{thm}
\begin{question}
It seems very probable that the conclusion of Theorem~\ref{thm1.3}
remains true if one replaces the hypothesis that the elements of
$\Sigma$ be special by the existence of just one special point on each
irreducible component of the Zariski closure of~$\Sigma$. The idea is
that the proof we give actually becomes easier if the Galois orbits
in~$\Sigma$ are bigger, and that is just what happens if instead of
special points we take non-special points.
\end{question}
We end this introduction with some words on the history of our proof,
on how it relates to the proof of more general cases, and on
perspectives of future research. In~\cite{Andre} André has proved
Theorem~\ref{thm1.2} \emph{unconditionally} in the case of a product
of \emph{two} modular curves. It may be possible to use his method to
give an unconditional proof of~Theorem.~\ref{thm1.2}. We have not
tried to do so, as our main goal was to test our approach to the
André--Oort conjecture for higher dimensional sub-varieties in at least
one situation. The proofs of Theorems~\ref{thm1.2} and~\ref{thm1.3}
were obtained in March 1999, but writing it all up has been delayed
for some time. One reason for that was that more important cases of
the Andr\'e--Oort conjecture have been dealt with first: Hilbert
modular surfaces in~\cite{Edixhoven2} and a result on curves in
arbitrary Shimura varieties in~\cite{EdixhovenYafaev1}. The importance
of the last result is its application to transcendence theory. Yafaev
has extended~\cite{EdixhovenYafaev1} to the case of arbitrary curves
in Shimura varieties, assuming~GRH (see~\cite{Yafaev2}), and he has
generalised a result of Moonen to arbitrary Shimura varieties
in~\cite{Yafaev3}. In the mean time, Breuer has succeeded in adapting
the arguments of this article to the case of Drinfel'd modular curves
in positive characteristic, see~\cite{Breuer1} and~\cite{Breuer2}.

We hope that the more or less explicit methods of this article can be
generalized and combined with the more abstract ones
of~\cite{EdixhovenYafaev1} in order to treat the general case of the
Andr\'e--Oort conjecture, i.e., higher dimensional cases in general
Shimura varieties. An interesting problem that suggests itself is to
generalize Proposition~\ref{prop4.2}, i.e., to get an effective
criterion for irreducibility for images under Hecke
correspondences. Is there an effective version of the theorem by Nori
that is used in~\cite{EdixhovenYafaev1}? Another important problem is
to find good enough lower bounds for Galois orbits. This is the main
subject of~\cite{Yafaev2}. Can one make use of reduction modulo $p$,
as in~\cite{Moonen2}? For relations between the Andr\'e--Oort
conjecture and equidistribution properties we refer
to~\cite{Ullmo1}. Particularly interesting is the main result
of~\cite{Clozel-Ullmo1} concerning the equidistribution of ``strongly
special'' sub-varieties: in our case the special curves in~$\CC^n$
that project surjectively to all coordinates are ``strongly special''.

\section{Determination of the sub-varieties of Hodge type.}
\label{sec2}
The definition of the notion ``sub-variety of Hodge type of a Shimura
variety'' that we use is that of Moonen, see
\cite[Def.~1.1]{Edixhoven2}, or \cite[6.2]{Moonen2} and
\cite[Prop.~2.8]{Moonen3}.
\begin{prop}\label{prop2.1}
The sub-varieties of Hodge type of a finite product of modular curves
are precisely the special sub-varieties as defined in
Definition~\ref{def1.1}.
\end{prop}
\begin{proof}
Let $X$ be a product of modular curves (indexed by some finite
set~$S$) as in Definition~\ref{def1.1}.  Let $G$ denote the algebraic
group $\PGL_{2,\QQ}$ over $\QQ$, and let $\HHpm$ denote the double
half plane $\PP^1(\CC)-\PP^1(\RR)$. Then $X$, with its modular
interpretation, is a component of the Shimura variety associated to
the Shimura datum $(G^S,(\HHpm)^S)$, together with a suitable compact
open subgroup $K$ of~$G^S(\AAf)$.

By definition, the sub-varieties of Hodge type of $X$ are given by
triplets $(H,Y,g)$, with $(H,Y)$ a sub Shimura datum of
$(G^S,(\HHpm)^S)$, and $g$ an element of~$G^S(\AAf)$. Here $H$ is a
reductive subgroup of $G^S$, and $Y$ is an $H(\RR)$-orbit in
$(\HHpm)^S$, consisting of $h\colon\SS\to G^S_\RR$ that factor through
$H_\RR$ (here $\SS$ is the real algebraic group $\CC^*$, and $\HHpm$
is to be viewed as the $G(\RR)$-conjugacy class of the morphism
$a+bi\mapsto(\begin{smallmatrix}a & -b\\b & a\end{smallmatrix})$ from
$\SS$ to~$G_\RR$). To be precise, the sub-varieties of Hodge type
associated to such a triplet $(H,Y,g)$ are the irreducible components
of the image of $Y$ under the map $Y\to X$, $y\mapsto\pi(y,g)$,
where $\pi$ is the quotient map from $(\HHpm)^S\times G^S(\AAf)$ to
$X\subset G^S(\QQ)\backslash((\HHpm)^S\times G^S(\AAf)/K)$.

Let $Z$ be a special sub-variety of~$X$. We want to show that it is of
Hodge type. Since products of sub-varieties of Hodge type are again of
Hodge type, we may assume that $Z$ is of one of the two forms as in
Definition~\ref{def1.1}. If $Z$ is a CM point $x$, one can take $H$ to
be the torus with $\QQ$-points $K_x^*/\QQ^*$, where $K_x$ is the
endomorphism algebra of an elliptic curve corresponding to~$x$ (the
choice of an element~$h\colon\SS\to G_\RR$ lying over~$x$ gives $H$ as
the smallest $\QQ$-subgroup of~$G$ through which $h$ factors). In the
second case, one can take $H$ to be $G$, embedded in $G^S$ by the
morphism that sends $g$ to the $g_sgg_s^{-1}$.

Suppose now that $Z$ is a sub-variety of Hodge type of~$X$. We want to
show that $Z$ is special. Let $(H,Y,g)$ be a triplet as above that
gives rise to~$Z$.  Since we are only interested in connected
components, we may and do assume $H$ to be connected (replace $H$ by
its connected component~$H^0$, note that the elements of~$Y$ factor
through~$H^0_\RR$, and replace $Y$ by one of the finitely many
$H^0(\RR)$-orbits of~$Y$). The connected reductive algebraic subgroups
of $G$ are $G$ itself, the trivial subgroup, and the one dimensional
tori. Hence the image of $H$ under any of the projections $p_s$ from
$G^S$ to $G$ is of one of these three kinds. We note that the trivial
subgroup cannot occur, because for any $h$ in $(\HHpm)^S$, the
morphism $p_sh$ from $\SS$ to $G_\RR$ is non-trivial. If $p_sH$ is all
of $G$, then $p_sY=\HH^\pm$, and if $p_sH$ is a one-dimensional torus,
then $p_sY$ is a point, because it is an orbit for the action of
$H(\RR)$ under conjugation; such a point is necessarily a CM point (it
is a sub-variety of Hodge type of dimension zero). Hence $S$ is the
disjoint union of $S'$ and $S''$, with $S'$ the set of $s$ with
$p_sH=G$. We have $Z=Z'\times Z''$, with $Z''$ a CM point, and with
$Z'$ a sub-variety of Hodge type in the product of the $X_{\Gamma_s}$
for $s$ in~$S'$, because images of a sub-variety of Hodge type under a
morphism of Shimura varieties induced by a morphism of Shimura data
are again of Hodge type. Hence we have reduced the problem of showing
that $Z$ is of Hodge type to the case where $p_sH=G$ for all~$s$.

Suppose that we have $p_{s,t}H\neq G^2$ for some pair $(s,t)$ with
$s\neq t$. Then Goursat's lemma (which says that the subgroups of a
product $A\times B$ are the inverse images of graphs of isomorphisms
from sub-quotients of $A$ to sub-quotients of $B$) implies that
$p_{s,t}H$ is $G$, embedded by the map $x\mapsto
(g_sxg_s^{-1},g_txg_t^{-1})$, for some $g_s$ and $g_t$ in~$G(\QQ)$.
It follows that $p_{s,t}Z$ is the one dimensional sub-variety of Hodge
type of $X_{\Gamma_s}\times X_{\Gamma_t}$ associated to the Shimura
datum $(G,\HHpm)$ with embedding into $(\HHpm)^2$ via
$\tau\mapsto(g_s\tau,g_t\tau)$. Hence $p_{s,t}Z$ is itself a modular
curve of the form $X_\Gamma$, embedded as a Hecke correspondence in
$X_{\Gamma_s}\times X_{\Gamma_t}$. So we can view $Z$ as a sub-variety
of Hodge type in the product of this $X_\Gamma$ and the $X_{\Gamma_u}$
with $u$ not in $\{s,t\}$. But then induction on the number of
elements of $S$ finishes the proof.

Finally suppose that $p_{s,t}H = G^2$ for all $s\neq t$. Then we have
$H=G^S$ (induction on the cardinality of~$S$, Goursat's lemma and the
fact that the normal subgroups of $G^S$ are the $G^T$ with $T$ a
subset of~$S$), and $Z=X$, hence $Z$ is special.
\end{proof}

\section{Some general principles.}\label{sec3}
In this section we list and prove some results that we use in the
proofs of Theorems~\ref{thm1.2} and~\ref{thm1.3}. 

We begin by giving a more intuitive description of the notion of
special sub-variety.  Definition~\ref{def1.1} implies that the special
sub-varieties of a product of two modular curves $X_1$ and $X_2$ are
the following: CM points $(x,y)$, fibers of a projection
$X_1\times\{y\}$ or $\{x\}\times X_2$ over a CM point, or the graph of
a Hecke correspondence between $X_1$ and $X_2$, or $X_1\times X_2$
itself. In particular, the special sub-varieties of $\CC^2$, with $\CC$
viewed as the $j$-line $\SL_2(\ZZ)\backslash\HH$, are the CM points,
fibers of a projection over a CM point, $\CC^2$ itself, or the image
in $\CC^2$ of the modular curve $Y_0(n)$ (parameterizing elliptic
curves with a cyclic subgroup of order $n$) for some $n\geq1$, under
the map sending $(E,G)$ to $(j(E),j(E/G))$.

Let us now look at the special sub-varieties of a product $X$ of any
number of modular curves $X_s=X_{\Gamma_s}$, in the notation
of~\ref{def1.1}. Let $Z$ be a special sub-variety of $X$, arising from
a partition $(S_1,\ldots,S_r)$ of~$S$. Let $S''$ be the subset of $S$
consisting of those $s$ such that $p_sZ$ is a CM point, and let $S'$
be its complement. Then $Z$ decomposes as a product $Z'\times Z''$,
with $Z''$ a CM point, and $Z'$ projecting dominantly (surjectively,
in fact) to all $X_s$ (with $s$ in $S'$ of course). Now consider
projections $p_T\colon Z'\to X_T:=\prod_{s\in T}X_s$ for two element
subsets $T$ of~$S'$. Then $p_TZ'$ is either all of $X_T$, or it is the
graph of a Hecke correspondence, depending on whether $T$ meets two or
only one of the~$S_i$. Obviously, $Z'$ is contained in the
intersection $Z'''$ of the $p_T^{-1}p_TZ'$, for $T$ ranging over the
two element subsets of~$S'$. If we take one element $s_i$ in each
$S_i$ contained in~$S'$, then the projections of both $Z'$ and $Z'''$
to the product of the~$X_{s_i}$ are finite and surjective. Hence $Z'$
and $Z'''$ have the same dimension, and $Z'$ is actually an
irreducible component of~$Z'''$. Let us state this conclusion in the
following proposition.

\begin{prop}\label{prop3.1} 
Let $n\geq0$ be an integer. A closed irreducible sub-variety $Z$ of a
product $X$ of $n$ modular curves $X_1,\ldots,X_n$ is special if and
only if (1) all images of $Z$ under projection to one or two factors
are special, and (2) $Z$ is an irreducible component of the
intersection of the inverse images of its images under these
projections. Equivalently, the special sub-varieties of $X$ are the
irreducible components of loci defined by conditions that demand
certain coordinates to be CM points, and by the existence of an
isogeny of a given degree between certain pairs of coordinates.
\end{prop}
The following two lemmas follow directly from this proposition. 

\begin{lem}\label{lem3.2}
Let $n\geq0$ be integer, and let $\Gamma_i$ and $\Gamma_i'$ be
congruence subgroups of $\SL_2(\ZZ)$ for $i$ in $\{1,\ldots,n\}$, such
that $\Gamma_i'$ is contained in $\Gamma_i$ for every~$i$. Let $X$ be
the product of the $X_{\Gamma_i}$, and $X'$ the product of the
$X_{\Gamma_i'}$. Let $\pi$ be the morphism from $X'$ to $X$ induced by
the inclusions of the $\Gamma_i'$ in the~$\Gamma_i$. Let $Z$ be a
closed irreducible sub-variety of~$X$. Then the following statements
are equivalent:
\begin{enumerate}
\item $Z$ is special; 
\item every irreducible component of $\pi^{-1}Z$ is special;
\item at least one irreducible component of $\pi^{-1}Z$ is special. 
\end{enumerate}
\end{lem}
\begin{lem}\label{lem3.3}
Let $n$, the $\Gamma_i$ and $X$ be as in the preceding
proposition. Let $Z_1$ and $Z_2$ be two special sub-varieties
of~$X$. Then all irreducible components of $Z_1\cap Z_2$ are special.
\end{lem}
The notion introduced in the next definition will allow us to reduce
the proof of our main results to the case where the $\Gamma_i$ are
just $\SL_2(\ZZ)$, and where the Zariski closure of~$\Sigma$ is a
hyper-surface all of whose projections to products of all but one of
the $X_i$ are dominant.
\begin{defi}\label{defi3.4}
Let $k$ be a field, $n\geq0$ an integer, $X_1,\ldots,X_n$ curves over
$k$ (i.e., $k$-schemes of finite type, everywhere of dimension one).
For $I$ a subset of $\{1,\ldots,n\}$, let $p_I$ be the projection from
$X:=X_1\times\cdots\times X_n$ to $X_I:=\prod_{i\in I}X_i$. Let $Z$ be
a closed irreducible sub-variety of~$X$. A subset $I$ of
$\{1,\ldots,n\}$ is said to be minimal for $Z$ if $\dim(p_IZ)< |I|$,
but $\dim(p_JZ)=|J|$ for all $J$ strictly contained in~$I$; in this
case, $p_I$ is called a minimal projection for~$Z$.
\end{defi}
\begin{lem}\label{lem3.5}
Notation as in Definition~\ref{defi3.4}. Then $Z$ is an irreducible
component of the intersection of the $p_I^{-1}\ol{p_IZ}$, with $I$
minimal for~$Z$. 
\end{lem}
\begin{proof}
First of all, note that the problem is only about closed subsets,
hence we may and do replace all schemes here by their reduced
sub-schemes. We replace each $X_i$ by an irreducible component of it
that contains the image of $Z$ under~$p_i$. Let $U_i$ be affine open
in $X_i$, such that $U_i$ meets~$p_iZ$. For each $i$, let $t_i$ be a
regular function on $U_i$ that is transcendental over~$k$. After
renumbering the $X_i$, the elements $p_1^*t_1,\ldots,p_d^*t_d$ form a
transcendence basis over $k$ of the function field of~$Z$. Then, for
every $j>d$, $p_j^*t_j$ is algebraic over the first~$d$, and we find a
minimal set $I_j$ consisting of $j$ and the $i\leq d$ that occur in
the minimal dependence relation of~$p_j^*t_j$. It follows that the
intersection of the $p_I^{-1}\ol{p_IZ}$, with $I$ ranging over
the~$I_j$, is the union of a $d$-dimensional closed part
containing~$Z$, and another closed part whose image in
$X_1\times\cdots\times X_d$ has dimension less than~$d$. Hence $Z$ is
an irreducible component of that intersection. The intersection of all
$p_I^{-1}\ol{p_IZ}$ contains $Z$ and is contained in the intersection
that we just considered, hence has $Z$ as an irreducible component.
%%  question: is this other part of dimension at most d????
%%  is it there at all????
%%  I see that it is there in the following example: take the subvariety 
%%  of \CC^4 with coordinates x,y,z,t defined by xz-y and xt-y and t-z. 
%%  Then if one projects to the x,y plane, one finds the intersection 
%%   given by xz-y and xt-y, which has a two-dimensional component over 
%%  the point x=y=0. 
\end{proof}
\begin{prop}\label{prop3.6}
Let $Z$ be an irreducible closed sub-variety of a product
$X=X_1\times\cdots\times X_n$ of complex modular curves. Then $Z$ is
special if and only if for every subset $I$ of $\{1,\ldots,n\}$ that
is minimal for $Z$ we have: $|I|\leq2$ and $\ol{p_IZ}$ is special.
\end{prop}
\begin{proof}
This follows immediately from Proposition~\ref{prop3.1}. 
\end{proof}

\section{Special sub-varieties and Hecke correspondences.}
For integers $m\geq1$ and $n\geq1$ we let $T_m$ be the Hecke
correspondence on $\CC^n$ that sends a point $(j(E_1),\ldots,j(E_n))$
to the formal sum of the $(j(E_1'),\ldots,j(E_n'))$ with each $E_i'$ a
quotient of $E_i$ by a cyclic subgroup of order~$m$. In other words,
$T_m$ is the correspondence induced by the sub-variety of $\CC^n\times
\CC^n$ consisting of the $(x,y)$ such that, for every~$i$, $x_i$ and
$y_i$ are $j$-invariants of elliptic curves related by an isogeny with
kernel isomorphic to~$\ZZ/m\ZZ$. The aim of this section is to prove
the following theorem.
\begin{thm}\label{thm4.1}
Let $n\geq0$ be integer. Let $Y$ be a closed algebraic sub-variety of
$\CC^n$ all of whose irreducible components contain a special point
and are of the same dimension, $d$, say. Suppose that $Y$ is contained
in $T_mY$ for some integer $m>1$ composed of prime numbers~$l$ greater
than~$3$ and the degrees of the projections from the
irreducible components of~$Y$ to sub-products $\CC^d$ of~$\CC^n$. Then
every irreducible component of $Y$ is special.
\end{thm}
Before proving this theorem, we will establish some ingredients for
it. The idea is of course to use Proposition~\ref{prop3.6}. The
following proposition will be used to show that the subsets $I$ of
$\{1,2,\ldots,n\}$ that are minimal for an irreducible component $Z$
of $Y$ consist of at most two elements.
\begin{prop}\label{prop4.2}
Let $n\geq3$ be integer. Let $Z$ be a closed irreducible hyper-surface
in $\CC^n$, and suppose that all projections $p_I$ from $Z$ to
$\CC^{n-1}$ are dominant. Then for every integer $m>1$ composed of
prime numbers $l>3$ such that $l>\deg(p_I)$ for all~$I$, the image
$T_mZ$ of $Z$ is irreducible.
\end{prop}
\begin{proof}
Let $m$ be as in the proposition; we write it as $m=l_1^{e_1}\cdots
l_r^{e_r}$ with the $l_i$ distinct prime numbers, and with the
$e_i>0$. Let $G_i:=\SL_2(\ZZ/l_i^{e_i}\ZZ)/\{1,-1\}$, and let
$G:=G_1\times\cdots G_r$. Let $X$ be the modular curve corresponding
to this quotient $G$ of~$\SL_2(\ZZ)$. This curve $X$ parametrizes
elliptic curves with, for each~$i$, a symplectic level $l_i^{e_i}$
structure given up to sign. The group $G$ acts faithfully on $X$ with
quotient~$\CC$. We let $G^n$ act on $X^n$ and denote the quotient map
to $\CC^n$ by~$\pi_n$. Since $T_mZ$ is an image of $\pi_n^{-1}Z$, it
suffices to show that $\pi_n^{-1}Z$ is irreducible.

Let $V$ be an irreducible component of $\pi_n^{-1}Z$, and let $H$ be
its stabilizer in $G^n$ (i.e., $H$ is the subgroup of $g$ in $G^n$
such that $gV=V$). It suffices now to show that $H=G^n$, since then
$V=\pi_n^{-1}Z$, as $ G^n$ acts transitively on the set of irreducible
components of~$\pi_n^{-1} Z$.  Lemma~\ref{lem4.3} below says that it
is enough to prove that all projections from $H$ to $G^{n-1}$ are
surjective.  By symmetry, it suffices to consider the projection on
the first $n-1$ factors.  We consider the two diagrams:
$$
\xymatrix{
P \ar@(lu,ru)[]^{G^{n-1}} \ar[r] \ar@{}[rd]|*+{\square} \ar[d] & 
X^{n-1} \ar@(lu,ru)[]^{G^{n-1}} \ar[d]^{\pi_{n-1}} \\
Z \ar[r]^{p_{n-1}} & \CC^{n-1}
}
\qquad
\xymatrix{
V \ar@(lu,ru)[]^H \ar@{>->}[r]  & \pi_n^{-1}Z \ar@(lu,ru)[]^{G^n}
\ar[d] \ar@{>->}[r] & X^n \ar@(lu,ru)[]^{G^n} \ar[r] \ar[d] 
& X^{n-1}\times\CC \ar[r] \ar[d] \ar@(lu,ru)[]^{G^{n-1}\times\{1\}}
\ar@{}[rd]|*+{\square}
& X^{n-1} \ar@(lu,ru)[]^{G^{n-1}} \ar[d] \\
  & Z \ar@{>->}[r]& \CC^n \ar[r] & \CC^n \ar[r]^{p_{n-1}} 
& \CC^{n-1}
}
$$
with $P$ the fibered product, $\pi_{n-1}$ the quotient for the action
of $G^{n-1}$, and $p_{n-1}$ the projection from $Z$ to the first $n-1$
factors. All four morphisms in the Cartesian square defining~$P$ are
generically finite, and dominant. As the most right square is
Cartesian, $P$ is the inverse image of~$Z$ under the morphism from
$X^{n-1}\times\CC$ to~$\CC^n$. The canonical morphism of $\pi_n^{-1}Z$
to~$P$ is the quotient for the action by~$\{1\}\times G$; let $\ol{V}$
be the image of $V$ in~$P$ under this morphism. The morphism from $P$
to $Z$ is the quotient by~$G^{n-1}$. The fact the morphism
$\pi^{-1}Z\to P$ is $G^{n-1}\times\{1\}$-equivariant implies that the
stabilizer in~$G^{n-1}$ of $\ol{V}$ is the image $\ol{H}$ of $H$
in~$G^{n-1}$. 

But now our hypothesis that $l_i$ is greater than the degree $d$ of
the projection from $Z$ to $\CC^{n-1}$ imply that $P$ is
irreducible. Indeed, the $G^{n-1}$-set $\Irr(P)$ of irreducible
components of~$P$ has at most~$d$ elements, because each element has
degree at least one over~$X^{n-1}$. On the other hand, $G^{n-1}$ has
no proper subgroups of index at most~$d$ because each factor $G_i$ is
generated by its $l_i$-subgroups.  Hence $\ol{V}=P$ and
$\ol{H}=G^{n-1}$, which is just what we had to prove.
\end{proof}
\begin{lem}\label{lem4.3}
Let $G=G_1\times\cdots\times G_r$ be as above. Let $H$ be a subgroup
of $G^n$ with $n\geq2$ such that $p_IH=G^2$ for all projections
$p_I\colon G^n\to G^2$. Then $H=G^n$.
\end{lem}
\begin{proof}
Induction on~$n$. We may and do assume that $n\geq3$. We view $H$ as a
subgroup of the product of $G$ by~$G^{n-1}$. Then, by the induction
hypothesis, $H$ projects surjectively to both factors. Let $H_1:=H\cap
G^{n-1}$ and $H_2:=H\cap G$ (these are the kernels of the two
projections); these are normal subgroups of $G^{n-1}$ and $G$,
respectively. Goursat's Lemma then says that $H$ is the inverse image
of the graph of an isomorphism between $G^{n-1}/H_1$ and~$G/H_2$. 

The normal subgroups of $G$ are the kernels of the reduction morphisms
from $G$ to the products $\prod_i\SL_2(\ZZ/l_i^{f_i}\ZZ)/\{1,-1\}$
with $f_i\leq e_i$. To prove this, one first notes that
$\SL_2(\FF_l)/\{1,-1\}$ is simple for any prime $l\geq5$
(see~\cite[VIII, Thm.~8.4]{Lang2}). Then one lets $V_{l,e}$ be the
kernel of the reduction from $\SL_2(\ZZ/l^e\ZZ)/\{1,-1\}$ to
$\SL_2(\ZZ/l^{e-1}\ZZ)/\{1,-1\}$. As a representation of
$\SL_2(\FF_l)$, $V_{l,e}$ is isomorphic to $\Sym^2(\FF_l^2)$, hence is
irreducible for $l>2$. Finally, one uses that the $l$-torsion of
$\SL_2(\ZZ/l^e\ZZ)$ is contained in $V_{l,e}$ for $l\geq5$. Similarly,
the normal subgroups of $G^{n-1}$ are products of normal subgroups
of~$G$.

Suppose that $H_2\neq G$. We take $i$ such that $H_2\cap G_i$ is not
equal to~$G_i$. Then for a unique $j$ with $2\leq j \leq n$ the
intersection of the factor $G_i$ in the $j$th factor $G$ in $G^{n-1}$
with $H_1$ is not equal to~$G_i$. It follows that the projection
$p_{\{1,j\}}$ from $H$ to $G^2$ is not surjective, contradicting the
hypotheses of the Lemma. Hence $H_2=G$, $H_1=G^{n-1}$, and $H=G^n$.
\end{proof}
\begin{lem}\label{lem4.4}
Let $m\geq 2$ be an integer, $n\geq0$, and $x$ in~$\CC^n$. Then the
$T_m$-orbit $\cup_{i\geq0}T_m^ix$ is dense in $\CC^n$ for the
Archimedean topology.
\end{lem}
\begin{proof}
Since $T_m$ is the product of the correspondence (also denoted) $T_m$
on each factor~$\CC$, the proof is reduced to the case $n=1$. The
inverse image of $\cup_{i\geq0}T_m^ix$ under $j\colon\HH\to\CC$ is an
orbit of the subgroup $H$ of $\GL_2(\ZZ[1/m])$ generated by
$\SL_2(\ZZ)$ and the
element~$(\begin{smallmatrix}m&0\\0&1\end{smallmatrix})$. The little
computation:
$$
(\begin{smallmatrix}m^{-k}&0\\0&1\end{smallmatrix})
(\begin{smallmatrix}1&\ZZ\\0&1\end{smallmatrix})
(\begin{smallmatrix}m^k&0\\0&1\end{smallmatrix}) = 
(\begin{smallmatrix}1&m^{-k}\ZZ\\0&1\end{smallmatrix})
$$
shows that $H\cap\SL_2(\RR)$ is dense in~$\SL_2(\RR)$. (In fact, $H$
contains $\SL_2(\ZZ[1/m])$.)
\end{proof}

\begin{proof} (Of Theorem~\ref{thm4.1}.) 
Let $n$, $Y$, $d$ and $m$ be as in the statement of the theorem, and
let $Z$ be an irreducible component of~$Y$. 

We consider the correspondence $T_{m,Y}'$ from $Y$ to itself induced
by~$T_m$: if we view $T_m$ as a closed sub-variety
of~$\CC^n\times\CC^n$, then $T_{m,Y}'$ is given by the Cartesian
diagram:
$$
\xymatrix{
T_{m,Y}' \ar@{>->}[r] \ar@{}[rd]|*+{\square} \ar@{>->}[d] 
& Y\times Y \ar@{>->}[d] \\
T_m \ar@{>->}[r] & \CC^n\times \CC^n.
}
$$
The two projections from $T_{m,Y}'$ to~$Y$ are finite, because the
projections from $T_m$ to $\CC^n$ are finite.  As $Y$ is contained
in~$T_mY$, the two projections from $T_{m,Y}'$ to $Y$ are
surjective. We let $T_{m,Y}$ be the union of the irreducible
components of dimension~$d$ of~$T_{m,Y}'$. Then both projections from
$T_{m,Y}$ to $Y$ are finite and surjective.

Let $\Irr(Y)$ be the set of irreducible components of~$Y$. Then
$T_{m,Y}$ induces a correspondence $T_m$ on~$\Irr(Y)$. We replace $Y$
by the union of the irreducible components of it that lie in the
$T_m$-orbit of the element~$Z$ of~$\Irr(Y)$.  If $Z'$ is an
irreducible component of~$T_m Z$, then for any~$I$, $\dim p_I Z = \dim
p_I Z'$. Thus a subset $I$ of $\{1,2,\ldots,n\}$ is minimal for~$Z$ if
and only if it is minimal for all irreducible components of~$Y$.

Let $I$ be a subset of $\{1,2,\ldots,n\}$ that is minimal for~$Z$,
hence for all irreducible components $Z'$ of~$Y$. We want to apply
Proposition~\ref{prop4.2} to the~$\ol{p_IZ'}$, the closure in~$\CC^I$
of~$p_IZ'$, for all~$Z'$. For every $i$ in~$I$, the degree of the
projection from $\ol{p_IZ'}$ to~$\CC^{I-\{i\}}$ is at most that of the
projection from~$Z'$ to $\CC^J$ where $J\supset I-\{i\}$ is such that
$|J|=d$ and $\dim(p_JZ')=d$. Hence the hypotheses on~$m$ in
Proposition~\ref{prop4.2} are satisfied.

Suppose that $|I|\geq3$. Proposition~\ref{prop4.2} shows that all
$T_m\ol{p_IZ'}$ are irreducible. But then $\ol{p_IY}$ and
$T_m\ol{p_IY}$ have the same number of irreducible components, all
these components are of the same dimension, and $\ol{p_IY}$ is
contained in~$T_m\ol{p_IY}$. It follows that
$\ol{p_IY}=T_m\ol{p_IY}$. This contradicts the density of all
$T_m$-orbits in~$\CC^I$ (Lemma~\ref{lem4.4}). Hence we have
$|I|\leq2$.

In order to prove that $Z$ is special, it suffices to show that
$\ol{p_IZ}$ is special (Proposition~\ref{prop3.6}). If $|I|=1$ then
$p_IZ$ is a special point because $Z$ contains a special point.

Suppose now that $|I|=2$. We replace the $Y$ that we have
by~$\ol{p_IY}$, so our new $Y$ is a closed curve in $\CC^2$, with
quasi finite projections to both factors~$\CC$, of degrees $d_1$ and
$d_2$ that are less than each prime number $l$ dividing~$m$. We let
$T_{m,Y}$ be the correspondence from $Y$ to itself induced by~$T_m$ as
above. We would like to apply \cite[Theorem~6.1]{Edixhoven1}, but that
result only applies to irreducible curves in~$\CC^2$, and to $m$ that
are square free. We generalize the proof of that result to the present
situation. We start with \cite[Lemma~6.3]{Edixhoven1}.  Consider the
commutative diagram:
\begin{subeqn}
\xymatrix{
T_{m,Y} \ar[r]^{p_1} \ar[d] & Y \ar[d]^{p_1} \\
Y_0(m) \ar[r]^{p_1} & \CC
}
\end{subeqn}
All four maps in this diagram are quasi finite and dominant, and the
horizontal ones are finite and surjective. The hypotheses on the $l$
dividing $m$ imply that for each irreducible component $Z$ of~$Y$, the
fibered product of $Z$ and $Y_0(m)$ over $\CC$ is irreducible ($G$
does not have a proper subgroup of index at most $d_1$ or~$d_2$). It
follows that $T_{m,Y}$ maps surjectively to the fibered product of $Y$
and $Y_0(m)$ over~$\CC$. This means that for every $(x,y)$ on~$Y$,
$T_{m,Y}(x,y)$ surjects to~$T_m(x)$. For $n\geq1$ let $T_{m,Y,n}$ be
the correspondence on $Y$ obtained by taking in $T_{m,Y}^n$ the
irreducible components that correspond to isogenies with cyclic kernel
on the first coordinate, i.e., that send $(x,y)$ to the sum of those
$(x',y')$ in $T_{m,Y}^n(x,y)$ that correspond to isogenies $x\to x'$
with cyclic kernel (of order~$m^n$).

Let now $Z_0$ be in~$\Irr(Y)$. We can then choose elements~$Z_1$,
$Z_2$, etc. in~$\Irr(Y)$ such that $Z_i\subset T_{m,Y}Z_{i-1}$, such
that moreover $Z_i\subset T_{m,Y,i}Z_0$ (here we use that
$T_{m,Y}(x,y)$ surjects to~$T_m(x)$; when choosing the $Z_i$ we just
make sure that the isogeny on the first coordinate is cyclic). As
there are only finitely many possibilities for the~$Z_i$, it follows
that for some $n\geq1$ and some $Z$ in $\Irr(Y)$ we have:
$$
Z\subset T_{m,Y,n}Z.
$$
Let $T$ be the correspondence on $Z$ induced by~$T_{m,Y,n}Z$. Then, by
the same irreducibility argument as above for~$T_{m,Y}$, for each
$(x,y)$ in $Z$ the set $T(x,y)$ surjects to $T_{m^n}(x)$, with
$T_{m^n}$ the correspondence on $\CC$ given by isogenies with cyclic
kernel of order~$m^n$. By Lemma~\ref{lem4.4} all $T_{m^n}$-orbits in
$\CC$ are dense. It follows that all $T$-orbits in $Z$ are not
discrete, as their projection to $\CC$ is dense.

Now the rest of the proof of \cite[Theorem~6.1]{Edixhoven1} can be
applied almost without change. Let $X$ be an irreducible component of
the complex analytic variety $\pi^{-1}Z$, where $\pi\colon
\HH^2\to\CC^2$ is the quotient for the action of~$\SL_2(\ZZ)^2$. Let
$G_X$ be the stabilizer of $X$ in $G:=\SL_2(\RR)^2$. Then Lemmas~6.6,
6.7, 6.8~and~6.9 of \cite{Edixhoven1} can be applied to~$X$. (In the
proof of Lemma~6.9 we do not know the second coordinate of the
elements $g_{i,j}$, but that information was not used anyway.) Please
note the erratum at the end of this article for a correction to the
end of the proof of Theorem~6.1 of~\cite{Edixhoven1}. The proof of
Theorem~\ref{thm4.1} is now finished (we have shown that the
irreducible component $Z$ of $Y$ is special, but as $T_{m,Y}$ acts
transitively on~$\Irr(Y)$, all irreducible components are special).
\end{proof}

\section{Galois action.}\label{sec5}
We recall very briefly some facts about the action of
$\Gal(\Qbar/\QQ)$ on the set of $j$-invariants of elliptic curves with
complex multiplications, i.e., on the special points of~$\CC$. For
some more details, see~\cite[\S2]{Edixhoven1}. Let $E$ be an
elliptic curve over $\CC$ with complex multiplications by a quadratic
imaginary field~$K$. Then $\End(E)=O_{K,f}=\ZZ+fO_K$ for some unique
integer $f\geq1$ (the conductor of the order $O_{K,f}$ in the maximal
order~$O_K$). For each automorphism $\sigma$ of~$\CC$ we have
$\End(\sigma E)\cong O_{K,f}$. The set $S_{K,f}$ of isomorphism
classes of complex elliptic curves with endomorphism ring isomorphic
to $O_{K,f}$ is a $\Pic(O_{K,f})$-torsor, hence finite. It follows
that $\Aut(\CC)$ acts on $S_{K,f}$ via $\Gal(\Qbar/\QQ)$. The action
of $\Gal(\Qbar/K)$ is given by the morphism from
$\Gal(\Qbar/K)\to\Pic(O_{K,f})$ that is unramified outside~$f$ and
that sends the Frobenius element at a maximal ideal $m$ not
containing~$f$ to the inverse of the class~$[m]$ of~$m$
in~$\Pic(O_{K,f})$. This morphism is surjective, hence we have, by the
Brauer-Siegel theorem~\cite[Ch.~XVI]{Lang1}:
$$
|\Gal(\Qbar/K){\cdot}j(E)| = |\Pic(\End(E))| =
|\discr(\End(E))|^{1/2+o(1)}, \quad |\discr(\End(E))|\to\infty.
$$
Let $l$ be a prime number that is \emph{split} in~$\End(E)$, i.e., for
which $\FF_l\otimes\End(E)$ is isomorphic as a ring
to~$\FF_l\times\FF_l$. Let $m$ be one of the two ideals in $\End(E)$
of index~$l$. Then $E$ is a quotient of its Galois conjugate $[m]E$
via an isogeny of degree~$l$ (if $E\cong\CC/\Lambda$, then
$[m]E\cong\CC/m\Lambda$). It follows that we have the inclusion of
subsets of~$\CC$:
$$
\Gal(\Qbar/K){\cdot}j(E) \subset T_l (\Gal(\Qbar/K){\cdot}j(E)).
$$

\section{Existence of small split primes.}\label{sec6}
The effective Chebotarev theorem of Lagarias, Montgomery and Odlyzko,
assuming GRH, as stated in~\cite[Thm.~4]{Serre1} and the second remark
following that theorem, plus a simple computation (see Section~5
of~\cite{Edixhoven1}) give the following result.
\begin{prop}
For $M$ a finite Galois extension of $\QQ$, let $n_M$ be its
degree, $d_M=|\discr(O_M)|$ its absolute discriminant, and for $x$
in $\RR$, let $\pi_{M,1}(x)$ be the number of primes $p\leq x$ that
are unramified in $M$ and such that the Frobenius conjugacy class
$\Frob_p$ contains just the identity element of~$\Gal(M/\QQ)$.  Then
for $M$ a finite Galois extension of $\QQ$ for which GRH holds and for
$x$ sufficiently big (i.e., bigger than some absolute constant) such
that:
$$
x > 2(\log d_M)^2(\log(\log d_M))^2,
$$
one has:
$$
\pi_{M,1}(x)\geq \frac{x}{3n_M\log(x)}. 
$$
\end{prop}
We will apply this result in the following situation. Let $n\geq1$ be
an integer, and let $K_1,\ldots,K_n$ be quadratic sub-fields of~$\Qbar$
for which GRH holds. Let $M:=K_1\cdots K_n$ be the composite of
the~$K_i$. Then:
$$
d_M\leq|\discr(O_{K_1}\otimes\cdots\otimes O_{K_n})| =
(d_{K_1}\cdots d_{K_n})^{2^{n-1}}. 
$$
On the other hand, for each $i$ we have embeddings $K_i\to M$, hence:
$$
d_M = |\Norm_{K_i/\QQ}(\discr(O_M/O_{K_i}))|{\cdot}d_{K_i}^{[M:K_i]}
\geq d_{K_i}.
$$
The preceding two inequalities mean that, for our purposes, $\log d_M$
is of the same order of magnitude as the maximum of the~$\log d_{K_i}$.

For each~$i$, let $R_i$ be an order of~$K_i$. The number of primes
dividing the discriminant of~$R_i$ is of order~$o(\log(\discr
R_i))$ (indeed, if $P(n)$ denotes the number of primes dividing a
positive integer~$n$, then one has $P(n)\log P(n)\leq n$).
It follows that if $\max\{|\discr(R_i)|\;|\;1\leq i\leq n\}$
is bigger than some absolute constant, then there are primes $l$ split
in each~$R_i$, such that:
$$
l\leq (\log \max\{|\discr(R_i)|\;|\; 1\leq i\leq n\})^{2+o(1)},
$$
where the $o(1)$ does not depend on the fields~$K_i$.

\section{The case of a curve.}\label{sec7}
In this section we prove Theorems~\ref{thm1.2} and~\ref{thm1.3} for
the one-dimensional irreducible components of the Zariski closure of a
set of special points on a product of modular curves. By
Proposition~\ref{prop3.6} and Lemma~\ref{lem3.2} it suffices to
consider closed irreducible curves $Z$ in~$\CC^2$ that contain
infinitely many special points. Assume one of the following two
conditions: the generalized Riemann hypothesis is true for imaginary
quadratic number fields, or the special points can be taken in one
isogeny class. Then we will prove that $Z$ is special.

Even though Theorem~\ref{thm1.2} has been proved in \cite{Edixhoven1}
for curves in a product of two modular curves, we reprove it here in a
somewhat simplified way. Namely, it turns out that the first step of
the proof given in \cite{Edixhoven1} can be skipped, i.e., the
arguments of \cite[\S3]{Edixhoven1} are not needed. We also prove the
variant Theorem~\ref{thm1.3} in this section (this variant was not
treated in~\cite{Edixhoven1}). We should also mention that the next
section reproves the results of this section, but we think that this
section serves well as a kind of warming up exercise for the more
complicated arguments of the next section.

If one of the two projections from $Z$ to $\CC$ is not dominant, then
$Z$ is the inverse image under that projection of a special point,
hence special. So we assume that both projections are dominant. As $Z$
contains a dense set of points with coordinates in~$\Qbar$ (the
special points), $Z$ is defined over a finite extension of~$\QQ$, and
therefore has only finitely many Galois conjugates. Let $Z_\QQ$ be the
closed irreducible algebraic curve in~$\AA^2_\QQ$ that by base change
to $\CC$ gives the union of the finitely many Galois conjugates
of~$Z$. Let $d_1$ and $d_2$ be the degrees of the two projections
to~$\AA^1_\QQ$.

At this point, we proceed directly to the arguments
of~\cite[\S4]{Edixhoven1}. Let $x=(x_1,x_2)$ be a special point
in~$Z_\QQ(\CC)$. Let $l$ be a prime number that is split in both
$\End(x_i)$. Then we have (see Section~\ref{sec5}):
$$
\Gal(\Qbar/\QQ){\cdot}x \subset Z_\QQ(\CC)\cap T_lZ_\QQ(\CC), 
\quad |\Gal(\Qbar/\QQ){\cdot}x|\geq\max\{|\Pic(\End(x_i))|\;|\; 1\leq
i\leq 2\}.
$$
On the other hand, the intersection number in
$\PP^1_\QQ\times\PP^1_\QQ$ of the closures of $Z_\QQ$ and $T_lZ_\QQ$
is $2d_1d_2(l+1)^2$ (the bidegree of $T_lZ_\QQ$
is~$((l+1)d_1,(l+1)d_2)$). Hence:
$$
|Z_\QQ(\CC)\cap T_lZ_\QQ(\CC)| \leq 2d_1d_2(l+1)^2, \quad \text{if the
 intersection is finite.}
$$
\begin{lem}\label{lem7.1}
With the notation of this section, there exists a special point
$x=(x_1,x_2)$ in $Z_\QQ(\CC)$ and a prime number~$l$ such that:
\begin{enumerate}
\item $l>\max\{3,d_1,d_2\}$;
\item $l$ splits in $\End(x_i)$ for both~$i$;
\item $2d_1d_2(l+1)^2<\max\{|\Pic(\End(x_i))|\;|\; 1\leq i\leq 2\}$.
\end{enumerate}
\end{lem}
\begin{proof}
Let $\Sigma$ be the set of special points in~$Z_\QQ(\CC)$. The
function $\Sigma\to\ZZ$ sending $x$ to $\max\{|\discr(\End(x_i))|\;|\;
1\leq i\leq 2\}$ is not bounded, because for each possible value for
the discriminant there are only finitely many elliptic curves. We
recall that we have assumed that either GRH holds for imaginary
quadratic fields, or that $Z$ contains infinitely many special points
in one isogeny class. Let us first deal with the second case. Then we
have two imaginary quadratic fields $K_1$ and $K_2$ (possibly the
same) and an infinite set $\Sigma$ of $x$ in~$Z_\QQ(\CC)$ with
$\End(x_i)$ an order in~$K_i$. Then
$|\Pic(\End(x_i))|=|\discr(\End(x_i))|^{1/2+o(1)}$ by a simple
argument. The classical Chebotarev theorem (see for example
\cite[Ch.~VIII, \S4]{Lang1}) asserts that the set of primes $l$ that
are split in $M=K_1K_2$ has natural density~$1/n_M$ (actually,
Dirichlet density is good enough here). We note that the number of
primes $l$ that divide $\discr(\End(x_i))$ is at most
$\log_2|\discr(\End(x_i))|$. Hence there do exist $x$ and $l$ as
claimed.

Let us now assume that GRH holds for imaginary quadratic fields. Then
we use the Brauer-Siegel theorem (see Section~\ref{sec5}), and the
application of the effective Chebotarev theorem from
Section~\ref{sec6}. 
\end{proof}
Let now $x$ and $l$ be as in Lemma~\ref{lem7.1}. Then the intersection
$Z_\QQ\cap T_lZ_\QQ$ cannot be finite. As $Z_\QQ$ is irreducible, we
have:
$$
Z_\QQ \subset T_lZ_\QQ.
$$
Theorem~\ref{thm4.1} now implies that all components of $(Z_\QQ)_\CC$
are special, hence in particular that $Z$ is special. This finishes
the proof of Theorems~\ref{thm1.2} and~\ref{thm1.3} in the case of a
curve.

\section{Producing special curves from special points.}\label{sec8}
We start the proof of Theorems~\ref{thm1.2} and~\ref{thm1.3} in the
case of sub-varieties of arbitrary dimension. This proof will also
reprove the case of curves that was treated in the previous section.
By Lemma~\ref{lem3.2} it suffices to consider closed irreducible
sub-varieties $Z$ of~$\CC^n$ of dimension~$d\geq1$, that contain a
dense set $\Sigma$ of special points.
\begin{thm}\label{thm8.1}
Let $Z$ be a closed irreducible sub-variety of dimension $d\geq1$
of~$\CC^n$. Assume that $Z$ contains a dense set $\Sigma$ of special
points, and that at least one of the following conditions holds: GRH
is true for imaginary quadratic fields, or $\Sigma$ can be taken to
lie in one isogeny class. Then for all but finitely many $x$
in~$\Sigma$, there is a special curve $C$ contained in $Z$ with $x$
in~$C(\CC)$.
\end{thm}
The curves $C$ will be obtained via repeated intersections of
sub-varieties with their image under a suitable Hecke correspondence,
until we get an inclusion as in Theorem~\ref{thm4.1}. In order to
control the degrees of the sub-varieties in question, we review some
facts on intersection theory before starting the proof. Appendix~A of
\cite{Hartshorne1} is a good reference for what we need. It may help
to note that we only need intersections with divisors, as
in~\cite{Faltings1} (see also~\cite{deJong1}).

Let $k$ be a field, and $n\geq 0$ an integer. Let
$\PP:=(\PP^1_k)^n$. As the Chow ring of~$\PP^1_k$ is $\ZZ[x]/(x^2)$,
with $x$ the class of a rational point, the Chow ring of $\PP$ is
$A:=\ZZ[\eps_1,\ldots,\eps_n]$, with $\eps_i^2=0$ for all~$i$, and
with $\ZZ$-basis the family of~$\eps_I=\prod_{i\in I}\eps_i$ indexed
by subsets $I$ of~$\{1,\ldots,n\}$. Let $Z$ be a closed $k$-irreducible
sub-variety of~$\PP$, and let $d$ be its dimension. We write its class
$[Z]$ in~$A$ as $\sum_Ia_I(Z)\eps_I$ (of course, for $|I|\neq n-d$ we
have $a_I(Z)=0$). The coefficient $a_I(Z)$ is the degree of the
projection $p_{\ol{I}}\colon Z\to (\PP^1_k)^{\ol{I}}$, with $\ol{I}$
the complement of~$I$ in~$\{1,\ldots,n\}$. We let
$a(Z):=\max_Ia_I(Z)$. Let us suppose that $x$ is a closed point of
$\PP$ that does not lie on~$Z$, and that $k$ is not finite. Then we
want to produce a hyper-surface~$H$ of~$\PP$ that contains~$Z$,
avoids~$x$, and has all $a_I(H)$ suitably bounded in terms of
the~$a_I(Z)$. (The exact bound does not matter much; what is important
is that the bound is polynomial in the~$a_I(Z)$.) The line bundle
$\calL:=\calO(1,\cdots,1)$ on~$\PP$ is very ample, and its space of
sections gives an embedding of~$\PP$ in some projective
space~$\PP^N_k$. The class $[\calL]$ of $\calL$ in $A$ is
$\sum_i\eps_i$. The image of~$Z$ in $\PP^N_k$ has degree:
$$
[Z]{\cdot}[\calL]^d = [Z]{\cdot}(\eps_1+\cdots+\eps_n)^d = 
[Z]{\cdot}d!\sum_{|J|=d}\eps_J = d!{\cdot}\sum_{|I|=n-d}a_I(Z).
$$
In $\PP^N_k$ we can project $Z$ birationally onto a hyper-surface in
some $\PP^{d+1}_k$ of the same degree, that avoids the image
of~$x$. It follows that we can take $H$ such that:
$$
[H]=\left(d!{\cdot}\sum_{|I|=n-d}a_I(Z)\right)\sum_i\eps_i.
$$
For $l$ prime, let $[T_l]$ be the class in $A\otimes
A=\ZZ[\eps_1,\ldots,\eps_n,\eta_1,\ldots,\eta_n]$ of the
correspondence $T_l$ on~$\PP\times\PP$. As $T_l$ is the product of the
usual Hecke correspondences $T_l$ on each coordinate we get:
$$
[T_l] = (l+1)^n\prod_i(\eps_i+\eta_i).
$$
It follows that for $Z$ and $l$ as above we have:
$$
[T_lZ] = (l+1)^n[Z], \quad a_I(T_lZ) = (l+1)^na_I(Z)\quad\text{for
all~$I$}.
$$
\begin{proof}
(Of Theorem~\ref{thm8.1}.) As in the previous section, $Z$ has only
finitely many Galois conjugates, and we let $Z_\QQ$ be the closed
irreducible sub-variety of $\AA^n_\QQ$ such that $(Z_\QQ)_\CC$ is the
union of these Galois conjugates. Let $x=(x_1,\ldots,x_n)$ be
in~$\Sigma$, and let $l$ be a prime number that is split in each of
the~$\End(x_i)$. Then we have:
\begin{subeqn}\label{eqn8.1.1}
\Gal(\Qbar/\QQ){\cdot}x \subset Z_\QQ(\Qbar)\cap T_lZ_\QQ(\Qbar).
\end{subeqn}
Let $m_x:=\max\{|\discr(\End(x_i))|\;|\; 1\leq i\leq n\}$. As we only
need to prove a statement for all but finitely many~$x$, we may
suppose that $m_x$ is sufficiently large in terms of the~$a_I(Z_\QQ)$.
By the results of Section~\ref{sec6}, if we assume GRH, and by some
simple argument if $\Sigma$ lies in one isogeny class, we can take $l$
such that:
\begin{subeqn}\label{eqn8.1.2}
l<(\log m_x)^3,\quad l>3,  \quad\text{and $l>a_I(Z_\QQ)$ for all~$I$}.
\end{subeqn}
If $Z_\QQ$ is contained in $T_lZ_\QQ$ then $Z$ is special by
Theorem~\ref{thm4.1}, and Definition~\ref{def1.1} implies the
existence of a special curve $C$ as desired. Suppose now that $Z_\QQ$
is not contained in~$T_lZ_\QQ$. The results in Section~\ref{sec5} tell
us that, ignoring finitely many of the~$x$:
\begin{subeqn}\label{eqn8.1.3}
|\Gal(\Qbar/\QQ){\cdot}x| > m_x^{1/3}.
\end{subeqn}
The discrepancy between $\log m_x$ and $m_x$ will be heavily exploited
in the sense that, for $m_x$ large enough, any fixed power of $\log
m_x$ is less than~$m_x^{1/3}$. In particular, if $m_x$ is sufficiently
large with respect to $n$ and $a(Z)$, then the size of
$\Gal(\Qbar/\QQ){\cdot}x$ is larger than all intersection numbers that
we will encounter in the rest of the proof. 

Let $H$ be a hyper-surface in $(\PP^1_\QQ)^n$ that contains~$T_lZ_\QQ$,
that does not contain~$Z_\QQ$ and that satisfies:
\begin{subeqn}\label{eqn8.1.4}
[H]=\left((l+1)^nd!\sum_{|I|=n-d}a_I(Z)\right)\sum_i\eps_i.
\end{subeqn}
Let $Z_1$ be a $\QQ$-irreducible component of $Z_\QQ\cap H$ that
contains~$x$. We note that
$\dim(Z_1)=d-1$. Equations~(\ref{eqn8.1.1})--(\ref{eqn8.1.4}) imply
that $Z_\QQ\cap H$ is not finite: it contains more points, namely the
Galois orbit of~$x$, than the intersection number. Hence $d>1$ (so if
$d=1$ then we have proved that $Z$ is special).
%% A small computation shows:
%% \begin{subeqn}\label{eqn8.1.5}
%% a(Z_1) < d!\;n2^na(Z_\QQ)^2(\log m_x + 1)^{3n}.
%% \end{subeqn}

The idea is now to apply to $Z_1$ the same constructions as we have
just applied to~$Z_\QQ$. At this point we do not know whether $Z_1$
has a dense subset of special points, but we know that $Z_1$ contains
$\Gal(\Qbar/\QQ){\cdot}x$. We get a prime number~$l_1$ of size about
$a(Z_1)$ that is split in each of the~$\End(x_i)$. If $Z_1$ is
contained in $T_{l_1}Z_1$ then we get a special curve~$C$ in~$Z$ as
desired. If not, then we take a hyper-surface $H_1$ in $(\PP^1_\QQ)^n$
that contains $T_{l_1}Z_1$ but does not contain~$Z_1$ and that has
suitably bounded degree as above, and let $Z_2$ be an irreducible
component of $Z_1\cap H_1$, etc. As the dimension drops by one at each
intersection, we need to repeat this process at most $d-1$ steps. One
easily computes that $a(Z_i)$ is of order of magnitude at most the
$(3n)^i$th power of~$\log m_x$.  Hence the intersection $Z_i\cap H_i$
is never finite, which means that at some point we will have
$Z_i\subset T_{l_i}Z_i$, with $Z_i$ of dimension at least one. The
irreducible components of this $(Z_i)_\CC$ are then special
by~Thm.~\ref{thm4.1}, and we get a special curve $C$ as desired
in~$(Z_i)_\CC$.
\end{proof}

\section{End of the proof.}\label{sec9}
We will now finish the proof of Theorems~\ref{thm1.2}
and~\ref{thm1.3}. So let $Z$ be an irreducible component of the
Zariski closure of a set $\Sigma$ of special points in a finite
product of modular curves, and assume either GRH for imaginary
quadratic fields or that $\Sigma$ lies in one isogeny class. As
explained in the introduction, in the last case we may and do assume
that the isogenies preserve the product structure, i.e., the
$n$-tuples of elliptic curves corresponding to the elements
of~$\Sigma$ are isogeneous coordinate-wise.

We have to prove that $Z$ is special. Let $I$ be a subset of
$\{1,\ldots,n\}$ which is minimal for~$Z$ (see
Definition~\ref{defi3.4}). By Proposition~\ref{prop3.6} it suffices to
prove that $|I|\leq2$ and that $\ol{p_IZ}$ is special. If $|I|=1$ then
$\ol{p_IZ}$ is a special point, hence special. If $|I|=2$ then
$\ol{p_IZ}$ is a special curve as was proved in Section~\ref{sec7}
(and also in Section~\ref{sec8}).

So let us assume that $|I|\geq3$. We have to get a contradiction
now. We replace $Z$ by~$\ol{p_IZ}$, $n$ by~$|I|$, and renumber $I$ as
$\{1,2,\ldots,n\}$. By Lemma~\ref{lem3.2}, it suffices to consider the
case where the congruence subgroups are maximal, i.e., where $Z$ is
contained in~$\CC^n$. So now $n\geq3$, and $Z$ is an irreducible
hyper-surface in $\CC^n$ all of whose projections to coordinate
hyperplanes are dominant. Theorem~\ref{thm8.1} tells us that there is
a Zariski dense subset $\calY$ of special curves in~$Z$. For each $Y$
in $\calY$ we let $I_Y$ be the set of $i$ in $\{1,\ldots,n\}$ such
that the projection $p_i\colon Y\to\CC$ is surjective (note that as
the $Y$ in~$\calY$ are special, projecting surjectively or dominantly
under~$p_i$ is equivalent). The $I_Y$ are non-empty because each $Y$
is a curve. As there are only finitely many possibilities for~$I_Y$,
there is a subset $I$ of $\{1,\ldots,n\}$ such that the set of $Y$
with $I_Y=I$ is Zariski dense. We pick such a subset~$I$ and replace
$\calY$ by the set of $Y$ with $I_Y=I$. We renumber the set
$\{1,\ldots,n\}$ such that $I=\{1,\ldots,j\}$, with $j\geq1$.

We claim that $j\geq3$. Indeed, if $j=1$ then each $Y$ in~$\calY$ is
of the form $\CC\times\{x\}$ with $x$ in~$\CC^{n-1}$ special; the set
of these $x$ is then Zariski dense in~$\CC^{n-1}$, contradicting the
fact that $Z\neq\CC^n$. Let us suppose then that $j=2$. Then for each
$Y$ in~$\calY$ we have $Y = p_{1,2}Y\times p_{>2}Y$, with $p_{1,2}Y$ a
special curve in~$\CC^2$, and $p_{>2}Y$ a special point in~$\CC^{n-2}$
(here $p_{>2}$ denotes the projection on the last $n-2$
coordinates). The set of the special points~$p_{>2}Y$ is Zariski dense
in~$\CC^{n-2}$. But, over a non-empty Zariski open subset
of~$\CC^{n-2}$ the fibers of $Z$ under the projection to the last
$n-2$ coordinates are curves in~$\CC^2$, of a fixed degree. Hence,
after shrinking~$\calY$, we may assume that the degrees of the
$p_{1,2}Y$ are all equal. As the set of special curves in~$\CC^2$
that project surjectively under the two projections and that have a
fixed degree is finite, this contradicts the Zariski density of the
union of the~$p_{1,2}Y$. Hence we have $j\geq3$.

Let $x_1$ be any special point in~$\CC$. For simplicity we suppose
that the CM-field of~$x_1$ is different from $\QQ(i)$ and
$\QQ(\sqrt{-3})$. Let $Z'$ be the Zariski closure of the union of the
intersections $Y\cap(\{x_1\}\times\CC^{n-1})$, $Y$ ranging
over~$\calY$. We note that $Z'$ is contained in
$Z\cap(\{x_1\}\times\CC^{n-1})$, and that $Z'$ is the Zariski closure
of a set of special points, contained in one isogeny class if $\Sigma$
is contained in one isogeny class: the first $j$ coordinates of an
element of $Y\cap(\{x_1\}\times\CC^{n-1})$ are isogeneous to~$x_1$,
and $p_{>j}Y$ is a special point whose coordinates lie in fixed
isogeny classes. Hence, by induction on~$n$, all irreducible
components of~$Z'$ are special.

Let $Z'_1,\ldots,Z'_r$ be the irreducible components of~$Z'$.  For $Y$
in $\calY$ let $\wt{Y}\to Y$ denote the normalization morphism. After
suitably renumbering the $Z'_i$ a Zariski dense subset of the $Y$
in~$\calY$ have the property that the number of points on~$\wt{Y}$
mapping to $Y\cap Z'_1$ is at least $1/r$ times the number of points
on~$\wt{Y}$ mapping to~$Y\cap Z'$. We replace $\calY$ by such a
subset. As $Z'_1$ is contained in~$Z$, we have
$Z'_1\neq\{x_1\}\times\CC^{n-1}$ (recall that the dimension of~$Z$ is
$n-1$ and that $Z$ projects surjectively to all coordinate
hyperplanes). As $Z'_1$ is special and not equal to
$\{x_1\}\times\CC^{n-1}$, there are $i$ and~$j$ with $1<i<j$ such that
$p_{i,j}Z'_1$ is a strict special sub-variety $S$ of~$\CC^2$.  We renumber
the indices so that this is so for $i=2$ and $j=3$.

The $p_{\leq3}Y$, for $Y$ ranging through~$\calY$, form a set of
special curves in~$\CC^3$ with the property that under projections to
all coordinate hyperplanes they give a Zariski dense set of special
sub-varieties (curves or points) in~$\CC^2$. For each $Y$ in~$\calY$,
let $N_Y\colon \wt{p_{\leq3}Y}\to p_{\leq3}Y$ be the normalization
map. As both $\wt{Y}$ and $\wt{p_{\leq3}Y}$ are quotients of $\HH$ by
a congruence subgroup of~$\SL_2(\ZZ)$, the morphism $N_Y$ is finite
and locally free, and unramified at all points $z$ with
$p_1(N_Yz)=x_1$ (note that $x_1$ is different from~$0$ and~$1728$). It
follows that, for all $Y$ in~$\calY$, the number of $x$ on
$\wt{p_{\leq3}Y}$ with $p_1(N_Yx)=x_1$ and $p_{2,3}N_Yx\in S$ is at
least $1/r$ times the number of points $x$ on $\wt{p_{\leq3}Y}$ with
$p_1(N_Yx)=x_1$. The following lemma shows that this is not the case,
and the proof of Theorems~\ref{thm1.2} and~\ref{thm1.3} is finished.

\begin{lem}\label{lem9.1}
Let $\calY$ be a set of special curves in~$\CC^3$ that map
surjectively to~$\CC$ under projection to the first coordinate, and
dense in some sub-variety of~$\CC^3$ (possibly equal to~$\CC^3$) that
projects dominantly to all coordinate hyperplanes. For $Y$ in~$\calY$,
let $N_Y\colon \wt{Y}\to Y$ be the normalization map. Let $x_1$ be a
special point in~$\CC$ whose CM-field is not $\QQ(i)$
or~$\QQ(\sqrt{-3})$, let $S\subset\CC^2$ be a special point or curve,
and let $r\geq1$. Suppose that for each $Y$ in~$\calY$ the number of
$z$ in~$\wt{Y}$ with $p_1(N_Yz)=x_1$ and $p_{2,3}(N_Yz)\in S$ is at
least $1/r$ times the number of $z$ in~$\wt{Y}$ with
$p_1(N_Yz)=x_1$. Then one has a contradiction.
\end{lem}
\begin{proof}
We assume all hypotheses of the lemma, and will arrive at a
contradiction. For $Y$ in~$\calY$, we let $I_Y$ be as before. Then we
may suppose that there is a subset $I$ of~$\{1,2,3\}$ such that
$I_Y=I$ for all~$Y$. By our assumptions, $I$ contains~$1$. 

Let us suppose that $I=\{1\}$. Then $p_{2,3}Y$ is in~$S$ for each $Y$
in~$\calY$, contradicting the assumption that the union of the
$p_{2,3}Y$ is dense in~$\CC^2$. 

Let us then suppose that $I$ has cardinality~$2$, say $I=\{1,2\}$ (if
necessary, we renumber the last two coordinates). Then, for each $Y$
in~$\calY$, $p_3Y$ is a one element subset of~$\CC$, and
$p_{2,3}Y=\CC\times p_3Y$. It follows that $S$ is a special curve,
surjecting to~$\CC$ under the projection~$p_2$ (note that, with our
notation, we have $p_2\circ p_{2,3}=p_3$). Then, for all $Y$
in~$\calY$, the cardinality of $S\cap p_{2,3}Y$ is bounded by the
degree, $m$ say, of~$p_2\colon S\to\CC$. On the other hand, $\wt{Y}$
is of the form~$Y_0(n)$, mapped to~$\CC^3$ by sending the isomorphism
class of an isogeny $f\colon E_1\to E_2$ with $\ker(f)$ isomorphic to
$\ZZ/n\ZZ$ to $(j(E_1),j(E_2),x_3)$ with $\{x_3\}=p_3Y$. Let $E_1$ be
such that $j(E_1)=x_1$. Then there are $\psi(n):=|\PP^1(\ZZ/n\ZZ)|$
subgroups of $E_1$ that are isomorphic to~$\ZZ/n\ZZ$. Two such
subgroups of~$E_1$ give the same point on~$\wt{Y}$ if and only if they
are equal (recall that $\Aut(E_1)=\{1,-1\}$). Hence there are at least
$\psi(n)$ points on~$\wt{Y}$ with first coordinate~$x_1$. Suppose
now that two such subgroups both lead to isogenies $f_1$ and~$f_2$
from $E_1$ to the same~$E_2$. Let $E'$ be the quotient of $E_1$ by
$\ker(f_1)\cap\ker(f_2)$, and $f_1',f_2'\colon E'\to E_2$ the
resulting cyclic isogenies, of degree~$d$, say. Then $f_2'\circ
(f_1')^\vee$ is an endomorphism of~$E_2$, with kernel isomorphic
to~$\ZZ/d^2\ZZ$. This endomorphism, together with~$f_1$,
determines~$f_2$. The number of endomorphisms of $E_2$ with kernel
isomorphic to~$\ZZ/d^2\ZZ$ is at most~$2^{\pi(d)}$, where, for an
integer~$i$, $\pi(i)$ is the number of (distinct) prime numbers
dividing~$i$. As $d$ divides~$n$, we see that the number of $j(E_2)$
arising like this is at least $\psi(n)/2^{\pi(n)}$. It follows that
$\psi(n)/2^{\pi(n)}r \leq m$.  But then there are only finitely many
possibilities for~$n$, contradicting that the union of the $p_{1,2}Y$
is dense in~$\CC^2$.

Finally, let us suppose that $I=\{1,2,3\}$. Let $Y$ be
in~$\calY$. Then, for some integers $n_{1,2}$ and~$n_{1,3}$,
$p_{1,2}Y$ is the image of $Y_0(n_{1,2})$, and $p_{1,3}Y$ is the image
of $Y_0(n_{1,3})$. Considering the intersection of the kernels of the
corresponding isogenies of degrees $n_{1,2}$ and~$n_{1,3}$ shows that
there are unique positive integers $n_1$, $n_2$ and~$n_3$ such that
$\wt{Y}$ has the following moduli interpretation: $\wt{Y}$ is the set
of isomorphism classes of $(E,H_1,H_2,H_3)$ with $E$ a complex
elliptic curve, $H_i$ a subgroup of $E$ isomorphic to~$\ZZ/n_i\ZZ$,
such that for $i\neq j$ one has $H_i\cap H_j=\{0\}$. The map
$\wt{Y}\to\CC^3$ sends $(E,H_1,H_2,H_3)$ to the point with
coordinates~$j(E/H_i)$. In particular, we have $n_{i,j}=n_in_j$.  A
good way to see what happens here is to consider, for each prime
number~$p$, the tree of lattices in~$\QQ_p^2$ up to~$\QQ_p^*$, and to
note that three points in a tree define a unique ``center'': the point
from which the paths to the three given points have disjoint
edges. Also, one uses that the action of $\PGL_2(\ZZ_p)$ on the set of
infinite non self-intersecting paths from the class of the standard
lattice $[\ZZ_p^2]$ is $3$-transitive, in order to see that the
triplet $(n_1,n_2,n_3)$ determines~$Y$.

For each $(E,H_1,H_2)$ as above with $j(E/H_1)=x_1$ there is the same
number of possibilities for~$H_3$. The number of $p_2(N_Yz)$, with $z$
on~$\wt{Y}$ such that $p_1(N_Yz)=x_1$ is at least
$\psi(n_1n_2)/2^{\pi(n_1n_2)}$. As the set of~$n_1n_2$, for $Y$
varying in~$\calY$, is not finite, we may suppose that
$\psi(n_1n_2)/2^{\pi(n_1n_2)}>r$. It follows that $p_1S$ is not a
point. Similarly, $p_2S$ is not a point. Hence $S$ is a special curve
that projects surjectively to both factors, say with degree~$m$. We
will now show that this contradicts the fact that the set of~$n_2n_3$,
when $Y$ varies, is not bounded. Indeed, there exists an $(E,H_1,H_2)$
with $j(E/H_1)=x_1$, such that at least $1/r$ of the possibilities
for~$H_3$ give a $(j(E/H_2),j(E/H_3))$ in~$S$. This leads to at least
$\phi(n_3)/2^{\pi(n_3)}$ points $(j(E/H_2),j(E/H_3))$ in~$S$, with the
same first coordinate. Hence we have $\phi(n_3)/2^{\pi(n_3)}\leq r$,
and $n_3$ is bounded as $Y$ varies in~$\calY$. Similarly, $n_2$ is
bounded. But then $n_2n_3$ is bounded, contradicting the fact that the
union of the $p_{2,3}Y$ is dense.
\end{proof}

\begin{rem}\label{rem9.2}
The case $|I|=3$ in the proof of Theorems~\ref{thm1.2}
and~\ref{thm1.3} above admits a simpler argument. In that case, $Y$ is
a special curve in~$\CC^3$, hence of one of the two types treated in
the proof of Lemma~\ref{lem9.1}. As $Y$ is contained in~$Z$, the
projection of $Y$ to its image under a projection to~$\CC^2$ has a
degree that is bounded by the degree of the projection of $Z$
to~$\CC^2$. It follows that only finitely many $Y$ of the first type
(i.e., with all three projections surjective) are possible: if $Y$
corresponds to $(n_1,n_2,n_3)$, then the degree of the projection
from~$Y$ to its image under~$p_{i,j}$ is at least~$\phi(n_k)$, with
$\{1,2,3\}=\{i,j,k\}$. For $Y$ with one constant projection, say the
image of $Y_0(a)$ in~$\CC^2$, embedded in $\CC^3$ with third
coordinate~$x_3$, it follows that $a$ is bounded. This gives a
contradiction with the fact that the set $\calY$ is Zariski dense in
the hypersurface~$Z$.

It is not hard to generalize the description given in the proof of
Lemma~\ref{lem9.1} of special curves in $\CC^3$ that project
surjectively to $\CC$ under all (three) projections to the case of
curves in $\CC^n$ with that property in $\CC^n$ with $n$
arbitrary. One finds that the set of such curves is in bijection with
the set $\PGL_2(\QQ)\backslash(\PGL_2(\QQ)/\PGL_2(\ZZ))^n$, which one
can interpret as the set of relative positions of $n$ lattices
in~$\QQ^2$; the bijection is induced by the elements $g_i$ given in
Definition~\ref{def1.1}. For details see~\cite[\S1.3]{Breuer2}. It is
also interesting to observe that there are such curves that are not
contained in their image under Hecke correspondences $T_m$ of small
level (compared to their degree).
\end{rem}

\section{Erratum to~\cite{Edixhoven1}.}\label{sec10}
There are two minor things to be dealt with, both of which do not
invalidate the main result.

Serre has pointed out to me that it is used, in the proof of
\cite[Lemma~6.3]{Edixhoven1}, that for $p$ prime and at least~5,
$\SL_2(\FF_p)$ has no proper subgroup of index at most~$p$. This is
wrong, as for example $\SL_2(\FF_{11})$ contains a subgroup isomorphic
to $A_5$, which has index~$11$. But, as Galois wrote in his ``lettre
testament'', it is true for all $p>11$. Hence the 5 in Theorem~6.1
should be replaced with~$13$. Another way to fix this is to consider
only subgroups of index less than~$p$, as we do in this article. 

When refereeing Yafaev's thesis, Pink has observed that there is a gap
at the end of the proof of \cite[Theorem~6.1]{Edixhoven1}. In the
notation of that proof, there is an element $g$ in $\GL_2(\QQ)$ such
that the stabilizer in $\SL_2(\RR)\times\SL_2(\RR)$ of the irreducible
complex analytic sub-variety in $\HH\times\HH$ is the graph of the
automorphism of $\SL_2(\RR)$ given by conjugation by~$g$. The problem
is that in the last six lines of the proof it is assumed, without
justification, that $g$ has positive determinant. To repair this, we
replace the last twelve lines of the proof, i.e., starting at ``Let
$x$ be an element of~$X$.'', by what follows.
\begin{quote}
Let $x=(x_1,x_2)$ be an element of $X$ such that the two projections
from $X$ to $\HH$ induce isomorphisms on tangent spaces $T_X(x)\to
T_\HH(x_1)$ and $T_X(x)\to T_\HH(x_2)$. These tangent spaces are
naturally isomorphic to the quotients $\Lie(G_X)/\Lie(G_{X,x})$,
$\Lie(\SL_2(\RR))/\Lie(\SL_2(\RR)_{x_1})$ and
$\Lie(\SL_2(\RR))/\Lie(\SL_2(\RR)_{x_2})$. Since $X$ is a complex
analytic sub-variety of $\HH\times\HH$, the isomorphisms between the
tangent spaces are compatible with the complex structures. Write
$x_2=g'x_1$, with $g'$ in $\SL_2(\RR)$. Then $g'$ also induces an
isomorphism $T_\HH(x_1)\to T_\HH(x_2)$ of one-dimensional complex
vector spaces. It follows that conjugation by $g^{-1}g'$ on
$\SL_2(\RR)$ induces an automorphism of
$\Lie(\SL_2(\RR))/\Lie(\SL_2(\RR)_{x_1})$ that preserves
orientation. A simple computation shows that this implies that
$g^{-1}g'$ is in the connected component of identity of the normalizer
in $\GL_2(\RR)$ of $\SL_2(\RR)_{x_1}$, hence that $g$ has positive
determinant. (Note that $\SL_2(\RR)_{x_1}$ is a conjugate of
$\SO_2(\RR)$, whose normalizer in $\GL_2(\RR)$ is $\RR^*{\rm
O}_2(\RR)$, a group with exactly two connected components.) Hence $g$
has positive determinant, and we have $x_2=gx_1$. This means that
$X=\{(\tau,g\tau)\,|\,\tau\in\HH\}$. We may replace $g$ by multiples
$ag$ of it, with $a$ a non-zero rational number. So we can and do
suppose that $g\ZZ^2$ is contained in $\ZZ^2$ and that $\ZZ^2/g\ZZ^2$
is cyclic, say of order~$m$. It is now clear that $C$ is~$Y_0(m)$.
\end{quote}

\subsection*{Acknowledgments.}
It is a pleasure to thank Jean-Pierre Serre for the first correction
in the erratum above, and Richard Pink for his observation concerning
the second correction. I am also grateful to the referees: their
remarks and suggestions have lead to a considerable improvement of the
quality of this text.

\vfill
\noindent
Bas Edixhoven\\
Universiteit Leiden\\
Mathematisch Instituut\\
Postbus 9512\\
2300\ RA\ \ Leiden\\
The Netherlands

\medskip\noindent
edix@math.leidenuniv.nl

\end{document}